\documentclass[fleqn]{mat01}
\usepackage{times,mathtimy,amssymb}
\begin{document}

\setcounter{page}{203}
\firstpage{203}

\font\zz=msam10 at 10pt
\def\Box{\mbox{\zz{\char'244}}}

\newcommand{\f}{\frac}
\newcommand{\s}{\sigma}
\newcommand{\la}{\lambda}
\newcommand{\Si}{\Sigma}
\newcommand{\iy}{\infty}
\newcommand{\what}{\widehat}
\newcommand{\lgra}{\longrightarrow}

\newtheorem{theo}{Theorem}
\renewcommand\thetheo{\arabic{section}.\arabic{theo}}
\newtheorem{theor}[theo]{\bf Theorem}
\newtheorem{lem}[theo]{Lemma}
\newtheorem{pot}[theo]{Proof of Theorem}
\newtheorem{propo}{\rm PROPOSITION}
\newtheorem{rema}[theo]{Remark}
\newtheorem{defn}[theo]{\rm DEFINITION}
\newtheorem{exam}{Example}
\newtheorem{coro}[theo]{\rm COROLLARY}
\def\conjecture{\trivlist\item[\hskip\labelsep{\it Conjecture.}]}
\def\exammp{\trivlist\item[\hskip\labelsep{\it Example.}]}
\def\theoo{\trivlist\item[\hskip\labelsep{\bf Theorem.}]}

\newcommand{\R}{\mathbb R}%
\newcommand{\C}{\mathbb C}%
\newcommand{\D}{\mathbb D}%
\newcommand{\Z}{\mathbb Z}%
\newcommand{\Q}{\mbox{$\mathbb Q$}}%
\newcommand{\N}{\mathbb N}%

\newcommand{\cK }{\mbox{$ \mathcal K $}}%
\newcommand{\cL }{\mbox{$ \mathcal L $}}%
\newcommand{\cG }{\mbox{$ \mathcal G $}}%
\newcommand{\cN }{\mbox{$ \mathcal N $}}%
\newcommand{\cP }{\mbox{$ \mathcal P $}}%
\newcommand{\cR }{\mbox{$ \mathcal R $}}%
\newcommand{\sdp}{\mbox{\rule{. 1mm}{2mm}$\! \times$}}%

\newcommand{\Up }{\mbox{$ \upsilon $}}%
\newcommand{\aaa}{\mbox{$\{A_n\}^\infty_{n=1} \subseteq {\mathcal A}$}}%

\newcommand{\sss}{\mbox{$\displaystyle\sum^\infty_{n=1}$}}%
\newcommand{\ccc}{\mbox{$\displaystyle\coprod^\infty_{n=1}$}}%

\def\e{\mbox{\rm{e}}}
\def\d{\mbox{\rm{d}}}
\def\V{\mbox{\rm{\bf V}}}
\def\f{\mbox{\rm{\bf f}}}
\def\r{\mbox{\rm{\bf r}}}

\font\xxxxx=tir at 7.6pt
\def\ee{\mbox{\xxxxx{e}}}

\title{Unsteady Stokes equations: Some complete general solutions}

\markboth{A Venkatlaxmi, B S Padmavathi and T Amaranath}{Unsteady Stokes equations}

\author{A VENKATLAXMI, B S PADMAVATHI and T AMARANATH}

\address{Department of Mathematics and Statistics, University of Hyderabad,
Hyderabad~500~046, India\\
\noindent E-mail: tasm@uohyd.ernet.in}

\volume{114}

\mon{May}

\parts{2}

\Date{MS received 2 December 2003}

\begin{abstract}
The completeness of solutions of homogeneous as well as non-homogeneous
unsteady Stokes equations are examined. A necessary and sufficient
condition for a divergence-free vector to represent the velocity field
of a possible unsteady Stokes flow in the absence of body forces is
derived.
\end{abstract}


\maketitle

\section{Unsteady Stokes flows: Homogeneous case}

The equations governing the motion of an arbitrary unsteady Stokes flow
of an incompressible, viscous fluid in the absence of any body forces
are
\begin{align}
\rho \frac{\partial \V}{\partial t} &= -\nabla p + \mu \nabla^{2}
\V,\\[.2pc]
\nabla \cdot \V &= 0,
\end{align}
where $\V$ is the velocity, $p$ is the pressure, $\rho$ is the density and
$\mu$ is the coefficient of dynamic viscosity of the fluid. 
Equation~(1) can also be written as
\begin{equation}
\mu \left( \nabla^{2} - \frac{1}{\nu} \frac{\partial}{\partial t}
\right) \V = \nabla p,
\end{equation}
where $\nu = (\mu/\rho)$ is the kinematic coefficient of viscosity.
Taking divergence of eq.~(3) and making use of eq.~(2), it is easy to see
that the pressure is harmonic. Hence, on operating the Laplace operator
on eq.~(3), we find that the velocity vector satisfies the equation
\begin{equation}
\nabla^{2} \left( \nabla^{2} - \frac{1}{\nu} \frac{\partial}{\partial t}
\right) \V = 0.
\end{equation}

\subsection{\it A complete general solution of unsteady Stokes equations}

Let $(\V, p)$ be any solution of (2) and (3). We define
\begin{equation}
\psi = \int_{0}^{t} \d\tau \int\!\int\!\int_D p(\xi, \eta, \zeta, \tau)
\left( \frac{\exp [-r^{2}/(4 \nu(t - \tau))]}{{\nu}^{1/2}[4 (t -
\tau)]^{3/2}} \right) \d\xi\ \d\eta\ \d\zeta,
\end{equation}
where $D$ is the region of the flow, which is simply connected and
\begin{equation*}
r = [(x - \xi)^{2} + (y - \eta)^{2} + (z - \zeta)^{2}]^{1/2}.
\end{equation*}
Then
\begin{equation}
\left(\nabla^2 - \frac{1}{\nu} \frac{\partial}{\partial t}\right) \psi =
p(x, y, z, t).
\end{equation}
Hence
\begin{equation}
\nabla^2 \left( \nabla^2 - \frac{1}{\nu} \frac{\partial}{\partial t}
\right) \psi = 0.
\end{equation}
Substituting (6) in (3), we get
\begin{equation}
\left(\nabla^2 - \frac{1}{\nu} \frac{\partial}{\partial t} \right)
\left[\V - \frac{1}{\mu} \nabla \psi \right] = 0.
\end{equation}
Hence
\begin{equation}
\V = \Phi + \frac{1}{\mu} \nabla \psi,
\end{equation}
where
\begin{equation}
\left(\nabla^2 - \frac{1}{\nu} \frac{\partial}{\partial t} \right)\Phi =
0.
\end{equation}
Substituting the expression (9) in eq.~(2), we get
\begin{equation}
\nabla^2 \psi = -\mu \nabla \cdot \Phi.
\end{equation}
A general solution of eq.~(7) is of the form $\psi = \psi_1 + \psi_2$,
where
\begin{align}
&\nabla^2 \psi_1 = 0,\\[.2pc]
&\left(\nabla^2 - \frac{1}{\nu} \frac{\partial}{\partial t} \right)
\psi_2 = 0,
\end{align}
(please see appendix). 

Hence from eq.~(6), we have 
\begin{align}
p &= \left(\nabla^2 - \frac{1}{\nu} \frac{\partial}{\partial t} \right)
(\psi_1 + \psi_2)\nonumber \\[.2pc]
  &= - \frac{1}{\nu} \frac{\partial\psi_1}{\partial t}.
\end{align}
Using eqs~(11) and (12), we get
\begin{equation}
\nabla^2\psi_2 = -\mu\nabla\cdot\Phi,
\end{equation}
which implies that
\begin{equation}
\frac{1}{\nu}\frac{\partial\psi_2}{\partial t} =
-\mu\nabla\cdot\Phi.
\end{equation}
Therefore
\begin{equation}
\frac{\partial\psi_2}{\partial t} = -\mu\nu\nabla\cdot\Phi \quad
\mbox{or} \quad \psi_2 =-\int_0^t \mu\nu\nabla\cdot\Phi (x,y,z,s) \d s,
\end{equation}
as $\psi_2(x,y,z,0) = 0$. 

So
\begin{align}
\V &= \Phi + \frac{1}{\mu}\nabla\psi \nonumber \\[.2pc]
&=\Phi + \frac{1}{\mu}\nabla\psi_1 + \frac{1}{\mu}\nabla\psi_2\nonumber\\[.2pc]
&=\Phi + \frac{1} {\mu}\nabla\psi_1 - \nu\nabla\int_0^{t}
\nabla\cdot\Phi(x,y,z,s)\d s.
\end{align}
This solution is similar to the Naghdi--Hsu type of solution~\cite{NH} given
for steady Stokes equations.

\section{Non-homogeneous case}

Consider the equations of motion of an unsteady Stokes flow of a viscous,
incompressible fluid in the presence of a body force {\bf f} given by
\begin{align}
\rho \frac{\partial \V}{\partial t} &= -\nabla p+\mu
\nabla^2 \V + \f,\\[.2pc]
\nabla\cdot \V &= 0.
\end{align}
The function {\bf f} can be represented as~\cite{M}
\begin{equation}
\f = \nabla \chi+\nabla \times \nabla \times(\r P) +
\nabla \times (\r T),
\end{equation}
where $\chi,\ P$ and $T$ are the given scalar functions. Taking
divergence of eq.~(19) and making use of eqs~(20) and (21), we observe
that
\begin{equation}
\nabla^2(\chi - p)=0.
\end{equation}
Hence
\begin{equation}
p = p^{\prime} + \chi,\quad\mbox{where}\quad \nabla^2p^{\prime}=0.
\end{equation}
Let $S$ be defined as the region
\begin{equation*}
S = \{(r,\theta,\varphi)\!:r_1\leq r \leq
 r_2,0\leq \theta \leq \pi,0\leq \varphi \leq 2\pi\}.
\end{equation*}
Chadwick and Trowbridge~\cite{CT} have shown that if {\bf V} is a vector field
possessing partial derivatives of orders up to two which are H\"{o}lder
continuous on $S$ and $\nabla\cdot \V = 0$ on $S$, then there are
scalar functions $A$ and $B$ on $S$ such that
\begin{equation}
\V = \nabla \times\nabla \times(\r A) + \nabla \times(\r B),
\end{equation}
where $A$ and $B$ are solutions of
\begin{align}
LA &=-\r \cdot \V,\\[.2pc]
LB &=-\r \cdot (\nabla \times \V),
\end{align}
where $L$ is the transverse part of the Laplace operator except for the
factor $1/r^2$.

\begin{theoo}
{\it A complete general solution of eqs~$(19)$ and $(20)$ in $S$ is given
by
\begin{align}
\V &= \nabla \times\nabla \times (\r A) + \nabla \times (\r B),\nonumber\tag{24}\\[.2pc]
p &=p_{0}+\chi+\frac{\partial}{\partial r}\left\lbrace r\left[
P + \mu \left(\nabla^2-\frac{1}{\nu}\frac{\partial}{\partial t}\right)
A\right]\right\rbrace,
\end{align}
where $A$ and $B$ satisfy equations
\begin{align}
&\mu\nabla^2\left(\nabla^2-\frac{1}{\nu}\frac{\partial}{\partial t}\right)
A + \nabla^2P=0,\\[.2pc]
&\mu\left(\nabla^2-\frac{1}{\nu}\frac{\partial}{\partial t}\right)
B+T = 0,
\end{align}
respectively.}
\end{theoo}

\begin{proof}
Let $A$ be a solution of eq.~(25). Operating
$\nabla^2(\nabla^2-\frac{1}{\nu}\frac{\partial}{\partial t})$ on both
sides of (25) and making use of eq.~(20), we find that
\begin{align}
\nabla^2\left(\nabla^2-\frac{1}{\nu}\frac{\partial}{\partial t}\right)LA 
&= -\nabla^2\left(\nabla^2-\frac{1}{\nu}\frac{\partial}{\partial t}\right)
(\r\cdot \V)\nonumber\\[.2pc]
&= -\r \cdot \nabla^2\left(\nabla^2-\frac{1}{\nu}\frac{\partial}{\partial t}
\right)\V.
\end{align}
Here we used the identity that
\begin{equation*}
\nabla^{2}(\r \cdot \V) = 2\nabla\cdot \V + \r \cdot \nabla^2 \V.
\end{equation*}
Since the operators $\nabla^2, L$ and $\partial/\partial t$ commute, we
rewrite eq.~(30) as
\begin{equation}
L\nabla^{2} \left(\nabla^2-\frac{1}{\nu}\frac{\partial}{\partial t}\right)
A= -\r \cdot \nabla^2\left(\nabla^2-\frac{1}{\nu}\frac{\partial}{\partial t}
\right)\V.
\end{equation}
We rewrite {\bf f} as
\begin{equation}
\f = \nabla \chi+\nabla\left(P+r\frac{\partial P}{\partial r}\right) - 
\r \nabla^{2} P + \nabla \times (\r T).
\end{equation}
Taking $\nabla^2$ of eq.~(19) and making use of (22) and (32), we
get
\begin{align}
\hskip -4pc \mu \nabla^2\left(\nabla^2-\frac{1}{\nu}\frac{\partial}{\partial t}\right)
\V &= \nabla^2(\nabla p)-\nabla^2 \f \nonumber\\[.2pc]
\hskip -4pc &= -\nabla^2 \nabla \left(P + r\frac{\partial P}{\partial r}\right)+
\nabla^2(\r\nabla^{2} P) - \nabla^2 (\nabla \times (\r T))\nonumber\\ 
\hskip -4pc &= -\nabla \left[\frac{\partial}{\partial r}(r\nabla^2P)\right] + \r
\nabla^{4} P - \nabla \times (\r \nabla^{2} T).
\end{align}
Here, from eqs~(31) and (33), we have
\begin{equation*}
\mu L \nabla^2 \left(\nabla^2-\frac{1}{\nu}\frac{\partial}{\partial t}\right)
A =-L\nabla^2P,
\end{equation*}
or
\begin{equation}
L \left[\mu \nabla^2\left(\nabla^2-\frac{1}{\nu}\frac{\partial}{\partial t}
\right)A + \nabla^2P\right]=0.
\end{equation}
From \cite{CT}, this implies that $\mu
\nabla^2(\nabla^2-\frac{1}{\nu}\frac{\partial}{\partial t})A
+\nabla^2P=f(r),$ where $f$ is an arbitrary function of $r$. Without
loss of generality, $f(r)$ can be neglected as in the proof given in
\cite{PRA}. Then we observe that $A$ satisfies eq.~(28).

Similarly, by considering eq.~(26), we find that
\begin{equation*}
\mu L\left(\nabla^2-\frac{1}{\nu}\frac{\partial}{\partial t}\right)B=-LT,
\end{equation*}
or
\begin{equation}
L\left[\mu \left(\nabla^2-\frac{1}{\nu}\frac{\partial}{\partial t}\right)
B + T\right] = 0.
\end{equation}
Following the result given in \cite{CT}, it follows that
\begin{equation}
\left[\mu\left(\nabla^2-\frac{1}{\nu}\frac{\partial}{\partial t}\right)
B + T\right] = g(r),
\end{equation}
where $g$ is an arbitrary function of $r$. Without loss of generality,
$g(r)$ can also be neglected as before. Then we find that $B$ satisfies
eq.~(29).

Substituting (21) and (24) in eq.~(19), we find that
\begin{align*}
&\nabla \left[\chi-p+\frac{\partial}{\partial r}(r P)+\mu
\frac{\partial} {\partial r}(r\nabla^2A)-\rho\frac{\partial}{\partial r}
(rA_{t})\right]- \r[\nabla^2P\nonumber\\[.4pc]
&\quad\ +\mu\nabla^4A -\rho\nabla^2A_{t}] +\nabla\times [\r (T+\mu
\nabla^2 B-\rho B_{t})]=0.
\end{align*}
Since $A$ and $B$ satisfy eqs~(28) and (29), we have
\begin{equation*}
p = p_{0} + \chi+\frac{\partial}{\partial
r}\left\lbrace r\left[P+\mu\left(\nabla^2-\frac{1}{\nu}\frac{\partial}
{\partial t}\right)A\right]\right\rbrace,\tag{27}
\end{equation*}
where $p_{0}$ is a constant. Hence a complete general solution of
eqs~(19) and (20) is given by {\bf V} and $p$ as given in (24) and (27)
respectively, where $A$ and $B$ satisfy eqs~(28) and (29) respectively.
\hfill $\Box$
\end{proof}

\subsection{\it Alternative proof of completeness}

Let $(\V, p)$ be a solution of (19) and (20). Since $p-\chi$ is
harmonic, it is possible to find a harmonic function $\phi$ where
\begin{equation}
\phi+r\frac{\partial \phi}{\partial r}=p-\chi ,
\end{equation}
and a function $A_{p}$ such that
\begin{equation}
\mu \left(\nabla^2-\frac{1}{\nu}\frac{\partial}{\partial t}\right)
A_{p} = \phi-P.
\end{equation}
Such a function $A_{p}$ satisfies
\begin{equation}
\mu\nabla^2\left(\nabla^2-\frac{1}{\nu}\frac{\partial}{\partial t}\right)
A_{p} = -\nabla^2 P.
\end{equation}
We observe that
\begin{align}
\V_1 &= \nabla \times \nabla \times (\r A_p),\\[.2pc]
p &= \chi+\frac{\partial}{\partial r}\left\lbrace r\left[\mu\left(\nabla^2
-\frac{1}{\nu}\frac{\partial}{\partial t}\right)A_{p}+P\right]\right\rbrace,
\end{align}
is a particular solution of (19) and (20), where $\f = \nabla \times
\nabla \times (\r P)+\nabla\chi$. Hence every $p$ satisfying
$\nabla^2p-\nabla^2\chi=0$ represents pressure in (19) and (20) for which
a possible velocity field is given by (40). Therefore in the homogeneous
case, every harmonic $p$ represents the pressure of a possible Stokes
flow. Consider
\begin{equation}
\V_{2} = \V - \V_{1},
\end{equation}
where $\V_1$ is as given in (40). Then
\begin{align}
\left(\nabla^2-\frac{1}{\nu}\frac{\partial}{\partial t}\right)\V_{2}
&= \left(\nabla^2-\frac{1}{\nu}\frac{\partial}{\partial t}\right)\V
-\left(\nabla^2-\frac{1}{\nu}\frac{\partial}{\partial t}\right)\V_1\nonumber\\[.2pc]
&= 0,
\end{align}
as {\bf V} and $\V_1$ are both solutions of (19). Since
\begin{equation}
\nabla \cdot \V_2 = 0,
\end{equation}
we can write
\begin{equation}
\V_2= \nabla \times\nabla \times (\r A_{c}) + \nabla \times (\r B),
\end{equation}
where
\begin{align}
&LA_{c} = -\r \cdot \V_{2},\\[.2pc]
&LB = -\r \cdot (\nabla \times \V_{2}).
\end{align}
Then, from (43) and (46)
\begin{equation*}
\left(\nabla^2-\frac{1}{\nu}\frac{\partial}{\partial t}\right)LA_{c}= -\r 
\cdot\left(\nabla^2-\frac{1}{\nu}\frac{\partial}{\partial t}\right)\V_2 = 0,
\end{equation*}
or
\begin{equation}
L\left(\nabla^2-\frac{1}{\nu}\frac{\partial}{\partial t}\right)A_{c} = 0.
\end{equation}
Similarly 
\begin{align}
\mu L\left(\nabla^2-\frac{1}{\nu}\frac{\partial}{\partial t}\right)B &= \mu
\left(\nabla^2-\frac{1}{\nu}\frac{\partial}{\partial t}\right)LB \nonumber\\[.2pc]
&= -\mu \left(\nabla^2-\frac{1}{\nu}\frac{\partial}{\partial t}\right)({\bf
r}\cdot (\nabla \times \V_2))\nonumber\\ 
&=-LT.
\end{align}
Hence
\begin{equation}
\left(\nabla^2-\frac{1}{\nu}\frac{\partial}{\partial t}\right)A_{c} = 0
\end{equation}
and
\begin{equation*}
\mu\left(\nabla^2-\frac{1}{\nu}\frac{\partial}{\partial t}\right)B + T =
0.\tag{29}
\end{equation*}
If we write
\begin{align}
&A=A_{p}+A_{c},\\[.2pc]
&\mu\nabla^2\left(\nabla^2-\frac{1}{\nu}\frac{\partial}{\partial t}\right)A =
\mu\nabla^2\left(\nabla^2-\frac{1}{\nu}\frac{\partial}{\partial t}\right)
[A_{p}+A_{c}] = -\nabla^2P,
\end{align}
or
\begin{equation*}
\mu\nabla^2\left(\nabla^2-\frac{1}{\nu}\frac{\partial}{\partial t}\right)
A+\nabla^2 P=0.\tag{28}
\end{equation*}
Hence
\begin{align*}
\V &= \V_1 + \V_{2}\\[.2pc]
   &= \nabla \times \nabla \times (\r A_{p}) + \nabla \times \nabla 
\times(\r A_{c}) + \nabla \times(\r B) \\[.2pc]
&= \nabla \times \nabla \times(\r A) + \nabla \times(\r B)\tag{24}
\end{align*}
and
\setcounter{equation}{26}
\begin{equation}
p =p_{0}+\chi+\frac{\partial}{\partial r} \left\lbrace r \left[P + \mu\left(
\nabla^2-\frac{1}{\nu}\frac{\partial}{\partial t}\right)A\right]\right\rbrace,
\end{equation}
where
\begin{align}
&\mu\nabla^2\left(\nabla^2-\frac{1}{\nu}\frac{\partial}{\partial t}\right)
A + \nabla^2 P =0,\\[.2pc]
&\mu \left(\nabla^2-\frac{1}{\nu}\frac{\partial}{\partial t}\right)
B + T = 0.
\end{align}

\section{Condition for a possible Stokes flow}

We now derive a necessary and sufficient condition for a divergence free
{\bf V} satisfying
$\nabla^2(\nabla^2-\frac{1}{\nu}\frac{\partial}{\partial t})\V = 0$, to
be the solution of homogeneous, unsteady Stokes equations~(1) and (2).
We make the following observations:
\begin{enumerate}
\renewcommand{\labelenumi}{\arabic{enumi}.}
\leftskip -.2pc
\item On operating curl on either sides of eq.~(3), we find that
the vorticity $\nabla \times \V$ satisfies the equation
\setcounter{equation}{52}
\begin{equation}
\hskip -.55cm \left(\nabla^2-\frac{1}{\nu}\frac{\partial}{\partial t}\right)
(\nabla \times \V)= 0.
\end{equation}
\item If $(\V_{1}, p)$ and $(\V_{2}, p)$ are solutions of the
homogeneous, unsteady Stokes equations~(1) and (2), then ${\bf U} = \
\V_1 - \V_2$ satisfies the equation
\pagebreak
\begin{equation}
\hskip -.55cm \left(\nabla^2-\frac{1}{\nu}\frac{\partial}{\partial
t}\right){\bf U} = 0, 
\end{equation}
and $({\bf U}, p_0)$ is also a solution, where $p_0$ is a constant. 
\item If $(\V, p_1)$ and $(\V, p_2)$ are solutions of the
homogeneous, unsteady Stokes equations~(1) and (2), then $(p_1-p_2)$ is
constant. 
\item If $(\V, p)$ is a solution of (1) and (2) then $(\nabla
\times \V, 0)$ is also a solution.
\end{enumerate}

We now show that every divergence-free {\bf V} satisfying
\begin{equation*}
\nabla^2\left(\nabla^2-\frac{1}{\nu}\frac{\partial}{\partial t}\right)\V = 0,\tag{4}
\end{equation*}
need not be a solution of the homogeneous unsteady Stokes equations~(1)
and (2). Let us consider
\begin{equation}
\V = (y\hat{i} - x\hat{j})~\exp(\nu t),
\end{equation}
where $\hat{i}$ and $\hat{j}$ are the unit vectors in the cartesian 
coordinates $(x,y)$. It can be easily verified that $\nabla\cdot \V =
0$ and $\V$ satisfies eq.~(4).

Substituting {\bf V} in eq.~(3) and expressing the equations in a
component form we get
\begin{align}
p_x &= -\mu y~\exp(\nu t),\\[.2pc]
p_y &= \mu x~\exp(\nu t).
\end{align}
It can be seen that $p_{xy}\neq p_{yx}$. Therefore it is not
possible to find the pressure $p$ corresponding to the velocity
{\bf V} given in (55) and hence {\bf V} is not a solution of
the homogeneous, unsteady Stokes equations. We now derive the
necessary and sufficient condition for a divergence-free vector
{\bf V} satisfying eq.~(4) to be a possible velocity field
in an unsteady Stokes flow in the absence of any body forces.

\begin{theoo}
{\it Let \V\ satisfy
\begin{equation*}
\nabla^2\left(\nabla^2-\frac{1}{\nu}\frac{\partial}{\partial t}\right)\V = 0,\tag{4}
\end{equation*}
and
\begin{equation*}
\nabla \cdot \V = 0.\tag{2}
\end{equation*}
Then \V\ will represent the velocity of a possible unsteady Stokes flow
if and only if the vorticity $\nabla \times \V$ satisfies the equation
\begin{equation*}
\left(\nabla^2-\frac{1}{\nu}\frac{\partial}{\partial t}\right)
(\nabla \times \V) =0.\tag{53}
\end{equation*}}
\end{theoo}

\begin{proof}
Suppose {\bf V} is the velocity of a possible unsteady Stokes flow in
the absence of body forces. Then consider the Stokes equations
\begin{equation*}
\mu\left(\nabla^2 - \frac{1}{\nu}\frac{\partial}{\partial t}\right)\V 
= \nabla p.\tag{3}
\end{equation*}
We have already seen that on operating curl on both sides of the above
equation, we find that
\begin{equation*}
\left(\nabla^2 - \frac{1}{\nu}\frac{\partial}{\partial t}\right)
(\nabla \times \V) = 0.\tag{53}
\end{equation*}
Conversely, let $\nabla \times \V$ satisfy eq.~(53). Then we show that
{\bf V} represents a possible unsteady, Stokes flow. As {\bf V} is
divergence-free, following the result given in \cite{CT}, we can express {\bf
V} as
\setcounter{equation}{23}
\begin{equation}
\V = \nabla \times\nabla \times(\r A) + \nabla \times(\r B),
\end{equation}
where $A$ and $B$ satisfy the equations
\begin{align}
LA &=-\r \cdot \V,\\[.2pc]
LB &=-\r \cdot (\nabla \times \V).
\end{align}
Further, as {\bf V} and $\nabla \times \V$ satisfy eqs~(4) and
(53) respectively, the functions $A$ and $B$ will satisfy equations
$\nabla^2(\nabla^2-\frac{1}{\nu}\frac{\partial}{\partial t})A = 0$
and $(\nabla^2-\frac{1}{\nu}\frac{\partial}{\partial t})B = 0$ 
respectively. Then {\bf V} given in (24) and
\setcounter{equation}{57}
\begin{equation}
p = p_{0}+ \mu\frac{\partial}{\partial r} \left[r\left(\nabla^2-
\frac{1}{\nu}\frac{\partial}{\partial t}\right)A\right],
\end{equation}
is a solution of eqs~(1) and (2). Hence {\bf V} represents the velocity
of a possible unsteady Stokes flow in the absence of any body
forces.\hfill $\Box$ 
\end{proof}

In fact, as any solution of (4) can be expressed as $\V = \V_{1} + \V_{2}$ 
where
\begin{equation*}
\nabla^2 \V_1 = 0,\quad \left(\nabla^2 - \frac{1}{\nu}
\frac{\partial}{\partial t} \right)\V_{2} = 0,
\end{equation*}
we observe that from (53), the necessary and sufficient condition stated
above reduces to the condition that
\begin{equation*}
-\frac{1}{\nu}\frac{\partial}{\partial t}(\nabla \times \V_1) = 0.
\end{equation*}
Now we show that every harmonic $p$ represents the pressure in an
unsteady Stokes flow in the absence of any body forces.

For a given $p$, we can choose a function $\phi$ and a function $A_p$ as
in eqs~(37) and (38), where
\begin{equation}
\phi + r\frac{\partial\phi}{\partial r} = p
\end{equation}
and
\begin{equation}
\mu\left(\nabla^2-\frac{1}{\nu}\frac{\partial}{\partial t}\right)
A_p= \phi,
\end{equation}
where
\begin{equation*}
\nabla^2\left(\nabla^2-\frac{1}{\nu}\frac{\partial}{\partial t}\right)
A_p = 0.
\end{equation*}
Then we can show that
\begin{align*}
&\V_1 = \nabla \times \nabla \times({\bf r}A_p),\\[.2pc]
&p = \mu\frac{\partial}{\partial r} \left[r\left(\nabla^2-\frac{1}{\nu}
\frac{\partial}{\partial t}\right)A_p\right],
\end{align*}
is a possible unsteady Stokes flow in the absence of any body forces.

Observe that for a given velocity field {\bf V}, the pressure $p$ is
determined uniquely up to a constant. For a given harmonic $p$, $(\V_1,
p)$ is a particular solution of Stokes equations and $(\V, p)$ is a
solution of eqs~(1) and (2) if $(\V - \V_1)$ satisfies
\begin{equation*}
\left(\nabla^2-\frac{1}{\nu}\frac{\partial}{\partial t}\right)(\V - \V_1)= 0.
\end{equation*}

\section{Conclusions}

Some solutions of homogeneous and non-homogeneous, unsteady Stokes
equations have been discussed and their completeness has been
established. In the case of homogeneous equations, we found a solution
of the Naghdi-Hsu \cite{NH} type which is shown to be complete. However,
the solution that is derived in the more general case, i.e., in the case
of non-homogeneous equations in the presence of body forces, is more
advantageous as it does not involve an integral and yet is in a closed
form. It is also easy to employ this solution for boundary value
problems as the boundary conditions can be expressed more easily using
the scalar functions that appear in this solution. This fact has been
observed in the steady case also \cite{PRA,PNAU}. In fact, this solution
is also valid in an infinite domain which follows by employing the proof
given in \cite{PA}. Lastly, a necessary and sufficient condition has
been derived for a divergence-free vector field to be a solution of a
possible unsteady Stokes flow in the absence of body
forces.\vspace{-.5pc}

\section*{Appendix}

Suppose
\begin{equation*}
\nabla^2\left(\nabla^2- \frac{1}{\nu}\frac{\partial}{\partial t}\right)\psi = 0,\tag{7}
\end{equation*}
then
\begin{equation*}
\psi=\psi_{1}+\psi_{2},
\end{equation*}
where
\begin{align*}
&\nabla^2 \psi_{1}=0,\tag{12}\\[.2pc]
&\left(\nabla^2- \frac{1}{\nu}\frac{\partial}{\partial t}\right)\psi_{2}=0.\tag{13}
\end{align*}
\begin{proof}
Let
\renewcommand{\theequation}{A\arabic{equation}}
\setcounter{equation}{0}
\begin{equation}
\psi^{\prime}=\left(\nabla^2-\frac{1}{\nu}\frac{\partial}{\partial t}\right)\psi.
\end{equation}
This implies 
\begin{equation}
\nabla^2 \psi^{\prime}=0.
\end{equation}
Define
\begin{equation}
\psi_{2}=\psi+\int_0^{t}\nu \psi^{\prime}(x,y,z,s)\d s.
\end{equation}
Then
\begin{align*}
\left(\nabla^2- \frac{1}{\nu}\frac{\partial}{\partial t}\right)\psi_{2}
&= \left(\nabla^2- \frac{1}{\nu}\frac{\partial}{\partial t}\right)\psi
+\left(\nabla^{2} - \frac{1}{\nu}\frac{\partial}{\partial t}\right) 
\int_0^{t}\nu \psi^{\prime}(x,y,z,s)\d s\\[.2pc]
&= \psi^{\prime}-\psi^{\prime},\quad \mbox{using (A1) and (A2)}\\[.2pc]
&= 0.
\end{align*}
Let
\begin{equation}
\psi_{1}=-\int_0^{t}\nu \psi^{\prime}(x,y,z,s)\d s.
\end{equation}
Hence
\begin{equation*}
\psi = \psi_{1} + \psi_{2},
\end{equation*}
where
\begin{align*}
&\nabla^2 \psi_{1}=0,\tag{12}\\[.2pc]
&\left(\nabla^2- \frac{1}{\nu}\frac{\partial}{\partial t}\right)\psi_{2}=0.\tag{13}
\end{align*}
\end{proof}

\end{document}